\documentclass[final,3p,times]{elsarticle}

\usepackage{graphicx,amsmath,amssymb,txfonts}
\usepackage{algorithm}
\usepackage{algorithmic}
\usepackage{subfigure,color}
\usepackage[colorlinks]{hyperref}
\usepackage{hypernat}
\usepackage{multirow}
\usepackage{booktabs}
\usepackage{threeparttable}
\allowdisplaybreaks
\graphicspath{{Figs/}}

\newtheorem{theorem}{\textbf{Theorem}}

\newtheorem{lemma}{\textbf{Lemma}}

\newtheorem{example}{\textbf{Example}}
\graphicspath{{Figs/}}

\begin{document}

\begin{frontmatter}


\title{A linearized energy--conservative finite element method for
the nonlinear Schr\"{o}dinger equation with wave operator}


%
\author[HDU]{Wentao Cai}
\author[CUHK]{Dongdong He\corref{cor}}
\author[CSU]{Kejia Pan}
\address[HDU]{Department of Mathematics, School of Sciences, Hangzhou Dianzi University, Hangzhou 310018, China}
\address[CUHK]{School of Science and Engineering,
The Chinese University of Hong Kong, Shenzhen, 518172, China}
\cortext[cor]{Corresponding author. E-mail address:  hedongdong@cuhk.edu.cn}
\address[CSU]{School of Mathematics and Statistics, Central South University, Changsha 410083, China}

\begin{abstract}
In this paper, we propose a linearized finite element method (FEM) for solving the cubic nonlinear Schr\"{o}dinger equation with wave operator. In this method, a modified leap--frog scheme is applied for time discretization and a Galerkin finite element method is applied for spatial discretization. We prove that the proposed method keeps the energy conservation in the given discrete norm. Comparing with non-conservative schemes, our algorithm keeps higher stability. Meanwhile, an optimal error estimate for the proposed scheme is given by an error splitting technique. That is, we split the error into two parts, one from temporal discretization and the other from spatial discretization. First, by introducing a time--discrete system, we prove the uniform boundedness for the solution of this time--discrete system in some strong norms and obtain error estimates in temporal direction. With the help of the preliminary temporal estimates, we then prove the pointwise uniform boundedness of the finite element solution, and obtain the optimal $L^2$--norm error estimates in the sense that the time step size is not related to spatial mesh size. Finally, numerical examples are provided to validate the convergence-order, unconditional stability and energy conservation.
\end{abstract}

\end{frontmatter}
\section{Introduction}
\setcounter{equation}{0}

{In this paper, we study the following cubic nonlinear Schr\"{o}dinger equation with wave operator}
\begin{align}\label{Schrodinger}
u_{tt}-\Delta u+\mathbf{i}u_t+|u|^2u+w(\mathbf{x})u=0,
\end{align}
for $\mathbf{x} \in \Omega$ and $t \in [0,T]$, where $\Omega$ is a bounded
convex domain in $\mathbb{R}^d$ ($d=2,\ 3$), $\mathbf{i}=\sqrt{-1}$ is the complex unit, $w(\mathbf{x})$ is the real--valued  potential function. Meanwhile, the initial and boundary conditions are defined by
\begin{align}
& u(\mathbf{x},t)=0,
&& \textrm{$\mathbf{x} \in \partial \Omega$, $t \in [0,T]$,}
\label{BC}
\\[3pt]
& u(\mathbf{x},0)=u_0(\mathbf{x}), \quad u_t(\mathbf{x},0)=u_1(\mathbf{x}),
&& \textrm{$\mathbf{x} \in \Omega$},
\label{INIT}
\end{align}
where $u_0, u_1$ are two given functions.
An important property of equation (\ref{Schrodinger}) is the energy conservation.  Computing the inner product of Eq. (\ref{Schrodinger}) with $u_t$ in $L^2(\Omega)$, and taking the real parts, one can obtain the following energy conservative identity
\begin{align}
\frac{d}{dt}\int_{\Omega}\left(|u_t|^2+|\nabla u|^2+\frac{|u|^4}{2}+w(\mathbf{x})|u|^2\right)dx=0.
\end{align}

{The nonlinear Sch\"{o}dinger equation is one of the most important equations in mathematical physics, which is originated from
quantum mechanics. It has been widely used to model various nonlinear physical phenomena, such as underwater
acoustics~\cite{Tappert1977}, nonlinear optics~\cite{MalomB2005,NewellAC1985}, quantum condensates~\cite{HaseA} and other nonlinear phenomena~\cite{AbMJ1981}. The cubic nonlinear Sch\"{o}dinger equation is one of the most important models, and it is also known as Gross--Pitaevskii equation (GPE), which plays a fundamental role in modeling the hydrodynamics of Bose--Einstein condensate~\cite{Antoine2013,Bao2013,Erdos}.

Recently, nonlinear Schr\"{o}dinger--type equations have also been widely studied. The nonlinear Schr\"{o}dinger equation with wave operator is  one of most important nonlinear Schr\"{o}dinger--type equations, it has been derived from many physical areas.}  For example, the nonrelativistic limit of the Kelin--Gordon equation~\cite{Machihara,Schoene,Tsutumi}, the Langmuir wave envelope approximation in plasma~\cite{Berge} and the modulated planar pulse approximation of the sine--Gordon equation for light bullets~\cite{Bao,Xin}.

{ Due to the wide applications of the Schr\"{o}dinger and Schr\"{o}dinger--type equations, performing efficient and accurate numerical simulations  plays an essential role in many real applications.  It is remarkable that nonconservative schemes for Schr\"{o}dinger--type equations may lead to numerical blow--up~\cite{Zhang}. Therefore, conservative schemes become very important for Schr\"{o}dinger and Schr\"{o}dinger--type equations. In the last several decades, numerical simulations of both the nonlinear Schr\"{o}dinger equation and the nonlinear Schr\"{o}dinger--type equation have been studied extensively. For examples, finite difference methods~\cite{Zhang,He2017,Guo,ZhangLM02,WangTC06,ZhangLM03,ZhangLM12,Bao12,WangTC2015}, finite element methods \cite{JWang,WSun,JWang-1,CLC,CLC1,FCJ2008,Lu05,Xu05,Guo15} and Fourier spectral method~\cite{Wang07}. In the field of finite difference methods, an implicit nonconservative difference scheme had been developed in~\cite{Guo} for solving nonlinear Schr\"{o}dinger equation, the method needs lots of algebraic operators. Zhang et al.~\cite{Zhang} had pointed out that the nonconservative schemes may easily lead to the numerical solution blow up. Thus, a conservative difference scheme was provided for solving the nonlinear Schr\"{o}dinger in their work. Subsequently, conservative finite difference methods were developed to solve the Schr\"{o}dinger equation with wave operator~\cite{ZhangLM02,WangTC06,ZhangLM03,ZhangLM12,Bao12,WangTC2015}. In the field of finite element methods, Galerkin finite element methods were used to solve the generalized nonlinear Schr\"{o}dinger equation~\cite{JWang,CLC,CLC1}, the coupled nonlinear Schr\"{o}dinger equation~\cite{WSun}  and the nonlinear Schr\"{o}dinger--Helmholtz system~\cite{JWang-1}. In these works, optimal error estimates are achieved in the sense that the time step size is not related to spatial mesh size.   Moreover,  the local discontinuous Galerkin methods were used to simulate the 1D nonlinear Schr\"{o}dinger equation~\cite{Lu05,Xu05} and multi--dimensional nonlinear Schr\"{o}dinger equation with wave operator~\cite{Guo15}. However, in \cite{Guo15}, optimal error estimates were only obtained for the linear equation at semi--discrete level. Therefore, in the context of numerical analysis,  rigorous error estimates for the numerical scheme of nonlinear Schr\"{o}dinger equation with wave operator are needed.}



{Previously, when using Galerkin FEMs for solving partial differential equations (PDEs)\\
\cite{BB2000,Akri-Lar2005,EW1980}, to obtain the error estimates of linearized explicit (or semi--implicit), the pointwise uniform boundedness of numerical solution in certain strong--norms was often required. Traditionally, the inverse inequality and mathematical induction were applied to obtain the pointwise boundness of the numerical solution,
\begin{align*}
\|R_hu^n-U^n_h\|_{L^{\infty}}\le Ch^{-\frac{d}{2}}\|R_hu^n-U^n_h\|_{L^2}\le Ch^{-\frac{d}{2}}(\tau^p+h^{r+1}),
\end{align*}
where $u^n$ and $U^n_h$ are the exact and numerical solutions at time level $t_n$, respectively. $R_h$ is Ritz projection operator and $d$ is the dimension, $\tau, h$ are temporal and spatial mesh sizes, $p, r$ are positive integers, referring to the convergence order of the FEM. But, to get the uniform boundedness of numerical solutions, the above inequality results in an unnecessary restriction between the time step size and spatial mesh size. Recently, a new technique \cite{Li131,Li132} was introduced to analyze error estimates of the linearized semi--implicit FEMs for time--dependent nonlinear PDEs. In Li and Sun's works~\cite{Li131,Li132}, errors were split into two parts, one part was from the temporal discretization and the other part was from the spatial discretization. By analyzing the introduced time--discrete PDEs, the pointwise uniform boundedness of FEM solutions can be proved in the sense that there is no restriction between the time step size and spatial mesh size.  Comparing to previous error analysis with conditional stability in~\cite{BB2000,Akri-Lar2005,EW1980}, optimal error estimates were obtained unconditionally  in~\cite{Li131,Li132}. Recently, this new technique has been used to analyze linearized FEMs for nonlinear Schr\"{o}dinger type equation~\cite{JWang,WSun,JWang-1,CLC,CLC1} and many other PDEs~\cite{LiSun2014,LiWangSun2014,Gao2014,Gao20142,Ligaosun2014,LiZhang2017}. To our best knowledge, this new technique is mainly carried out for nonlinear parabolic type of equations and also possibly coupled with elliptic type of equations. However, when considering the nonlinear Schr\"{o}dinger equation with wave operator, it will possess the combination of both parabolic and hyperbolic properties,  and this could raise some complexity for the error analysis. }

In this paper, a linearized FEM is proposed to solve the cubic nonlinear Schr\"{o}dinger equation with wave operator (\ref{Schrodinger}) subject to initial boundary conditions (\ref{BC})--(\ref{INIT}). In this scheme, system \eqref{Schrodinger}--(\ref{INIT}) is discretized  by a modified frog--leap scheme in time direction and the Galerkin finite element method in spatial direction. The proposed scheme is a semi-implicit  linear method which  only needs to solve a linear system at each time step. Thus, the proposed scheme is simpler and more efficient than implicit nonlinear schemes, which need to do iteration at each time step.  More importantly, discrete energy of the proposed method is conserved so that the scheme will not yield blow up. Subsequently, we will apply the error splitting technique~\cite{Li131,Li132} to study the proposed linearized energy--conservative FEM. By introducing a time--discrete system, we  will prove the uniform boundedness of time--discrete solutions in certain strong norms, and give the error estimates of time--discrete solutions. Based on the uniform boundedness of time--discrete solutions and mathematical induction, we get the pointwise uniform boundedness of fully--discrete FEM solutions in the sense that there is no restriction between the time step size and spatial mesh size. With the help of the above  pointwise uniform  boundedness of the FEM solutions and the traditional error analysis method, we can obtain the optimal error estimates in $L^2$--norm. {Our work in this paper can be regarded as a complementary to Guo and Xu's work \cite{Guo15}}.

The rest of the paper is organized as follows. In section~\ref{sec2}, a linearized Galerkin FEM for the nonlinear Schr\"{o}dinger equation with wave operator  is given, together with some assumptions and notations. Meanwhile, the energy conservative property of the fully--discrete system is presented and the time--discrete system is introduced. In section~\ref{sec3}, the uniform boundedness of the time--discrete solutions is proved in some strong norms.    Moreover, the error estimates of the time--discrete solutions are obtained. In section~\ref{sec4}, based on the uniform boundedness of time--discrete solution in $H^2$--norm,  we obtain the uniform boundedness of the fully discrete solutions in $L^{\infty}$--norm. In section~\ref{sec5}, we derive the optimal error estimates in $L^2$--norm unconditionally. Section~\ref{sec6} shows the numerical results, which confirm well with the theoretical findings. Finally, a brief conclusion is provided in the final section.

\section{Linearized Galerkin FEM method}\label{sec2}

Let $\Omega$ be the convex polygon in $\mathbb{R}^2$ (or convex polyhedron in $\mathbb{R}^3$). Following the classical FEM theory, we define $\mathcal{T}_h$ be a quasi-uniform partition of $\Omega$ into triangles $\pi_h$ in $\mathbb{R}^2$ or tetrahedra in $\mathbb{R}^3$.
Let $h=\max\limits_{\pi_h\in\mathcal{T}_h}{\rm diam}\,\pi_h$ denotes the meshsize. Let $V_h$ be the finite--dimensional subspace of $H^1_0(\Omega)$, which consists of continuous piecewise polynomials of degree $r$ ($r\geq1$) on $\mathcal{T}_h$.

For any two complex functions $u, v$, the inner product is defined as
\begin{align}
(u,v)=\int_{\Omega}u(x)(v(x))^{*}dx,
\end{align}
where $v^*$ denotes the conjugate of $v$.

Let $R_h:H^1_0(\Omega)\rightarrow V_h$ be the Ritz projection operator defined by
\begin{align}\label{Ritz}
(\nabla (v-R_hv),\nabla \omega)=0, \;\forall\;\omega \in V_h.
\end{align}

By the classical FEM theory \cite{Thomee}, the following inequality is valid,
\begin{align}
&\| v-R_hv\|_{L^2}+h\|\nabla(v-R_hv)\|_{L^2}\leq Ch^{r+1}\| v\|_{H^{r+1}}, \label{theory}\\
&\|R_hv\|_{W^{1,p}}\leq C \|v\|_{W^{1,p}}, \quad p>1, \label{theory2}
\end{align}
for any $v\in H^{r+1}_0(\Omega)$.

In this paper, the following inverse inequality \cite{Thomee} is always used:
\begin{align}\label{inversein}
\| v\|_{L^{\infty}}\leq C h^{-\frac{d}{2}}\| v\|_{L^2},
\end{align}
for any $v\in V_h$ and $d=2,3.$

Let $\{t_n|t_n=n\tau; 0\leq n\leq N\}$ be a uniform partition of $[0,T]$ with time step size $\tau=\frac{T}{N}$, and $U^n= u(\cdot, t_n)$. We define
\begin{align*}
\delta^2_{\tau}U^n=\frac{U^{n+1}-2U^n+U^{n-1}}{\tau^2},\ D_{\tau}U^n=\frac{U^{n+1}-U^{n-1}}{2\tau}, \ U^{\overline{n}}=\frac{U^{n-1}+U^{n+1}}{2}.
\end{align*}

With above notations, a linearized leap--frog Galerkin FEM is to seek $U^{n+1}_h\in V_h$ such that
\begin{align}\label{FEMscheme}
(\delta^2_{\tau}U^n_h, v_h)+(\nabla U^{\overline{n}}_h,\nabla v_h)+\mathbf{i}(D_{\tau}U^n_h,v_h)+\left(|U_h^n|^2U^{\overline{n}}_h,v_h\right)+(wU^{\overline{n}}_h,v_h)=0, \ n=1,2,\cdots, N-1,
\end{align}
for any $v_h\in V_h.$ The initial and first step FEM solutions are defined by
\begin{align}\label{FEMIN}
U^0_h=R_h u_0,\ \ U^1_h=R_h\left(u_0+\tau u_1+\frac{\tau^2}{2}(\Delta u_0-\mathbf{i}u_1-|u_0|^2u_0-wu_0)\right).
\end{align}
Next, we give a theorem to show that {\color{blue}{Eqs. \eqref{FEMscheme}--\eqref{FEMIN}}} is an energy conservative scheme.

\begin{theorem}\label{Theorem1*}
The discrete energy of the FEM scheme Eqs. \emph{(\ref{FEMscheme})--(\ref{FEMIN})} is conservative, i.e.,
$$E^{N-1}=E^{N-2}=\cdots=E^1=E^0,$$
where the discrete energy is defined as
$$E^n\triangleq\left\|\frac{U^{n+1}_h-U^n_h}{\tau}\right\|^2_{L^2}+\frac{1}{2}(\|\nabla U^{n+1}_h\|^2_{L^2}+\|\nabla U^n_h\|^2_{L^2})+\frac{\| U^n_hU^{n+1}_h\|^2_{L^2}}{2}+\frac{(wU^{n}_h,U^{n}_h)+(wU^{n+1}_h,U^{n+1}_h)}{2}.$$
\end{theorem}
{\bf Proof.}\quad Putting $v_h=U^{n+1}_h-U^{n-1}_h$ in (\ref{FEMscheme}) and taking the real parts for the each term of the result equation, one has
\begin{align}\label{Re-equ-1}
&\textrm{Re}(\delta^2_{\tau}U^n_h, U^{n+1}_h-U^{n-1}_h)+\textrm{Re}(\nabla U^{\overline{n}}_h, \nabla(U^{n+1}_h-U^{n-1}_h))+\textrm{Re}({\bf{i}}D_{\tau}U^n_h, U^{n+1}_h-U^{n-1}_h)\nonumber\\
&+\textrm{Re}(|U_h^n|^2U^{\overline{n}}_h, U^{n+1}_h-U^{n-1}_h)+\textrm{Re}(wU^{\overline{n}}_h, U^{n+1}_h-U^{n-1}_h)=0.
\end{align}
If we denote the terms of left side of above equation by $I_i$ ($i=1,\cdots,5$), then it is easy to see that
\begin{align*}
&I_1=\left\|\frac{U^{n+1}_h-U^n_h}{\tau}\right\|^2_{L^2}-\left\|\frac{U^{n}_h-U^{n-1}_h}{\tau}\right\|^2_{L^2},\\
&I_2=\frac{1}{2}\left(\|\nabla U_h^{n+1}\|^2_{L^2}+\|\nabla U_h^{n}\|^2_{L^2}\right)-\frac{1}{2}\left(\|\nabla U_h^{n}\|^2_{L^2}+\|\nabla U_h^{n-1}\|^2_{L^2}\right),\\
&I_3=0,\\
&I_4=\frac{\|U^{n+1}_hU^n_h\|^2_{L^2}}{2}-\frac{\|U^{n}_hU^{n-1}_h\|^2_{L^2}}{2},\\
&I_5=\frac{(wU^{n+1}_h,U^{n+1}_h)+(wU^{n}_h,U^{n}_h)}{2}-\frac{(wU^{n}_h,U^{n}_h)+(wU^{n-1}_h,U^{n-1}_h)}{2}.
\end{align*}
Thus, from above results and \eqref{Re-equ-1}, we have
$$E^n=E^{n-1}, \quad n=1,\cdots,N-1.$$
The proof is complete. \quad $\Box$

%
%

In this paper, we assume that the solution to the initial boundary problem (\ref{Schrodinger})--(\ref{INIT}) exists and satisfies,
\begin{align}\label{assumption}
&\| u_0\|_{H^{r+3}}+\| (u_0)_{ttt}\|_{H^{2}}+\| u\|_{L^{\infty}([0,T]; H^{r+1})}+\|u\|_{L^{\infty}([0,T];W^{1,\infty})}+\| u_t\|_{L^{\infty}([0,T];H^{r+1})}\nonumber\\
&+\| u_{tt}\|_{L^{\infty}([0,T];H^{r+1})}+\|u_{tttt}\|_{L^{\infty}([0,T];H^1)}\leq M_0,
\end{align}
{where $r\ (r\geq1)$ is the order of the Galerkin FE space used in (\ref{FEMscheme}), $M_0$ is a positive constant depends only on $\Omega$. In addition, the potential function $w(\mathbf{x})$ is assumed to belong to $H^{r+1}(\Omega)$. Denote $M=\|u\|_{L^{\infty}([0,T];H^2)}$}, then by the above assumption, $M$ is a finite positive number which only depends on $\Omega$.

Now we introduce a corresponding time--discrete system for the Galerkin FEM scheme (\ref{FEMscheme}),
\begin{align}\label{time-discrete}
\delta^2_{\tau}U^n-\Delta U^{\overline{n}}+\mathbf{i}D_{\tau}U^n+|U^n|^2U^{\overline{n}}+wU^{\overline{n}}=0, \ n=1,2,\cdots, N-1,
\end{align}
with the following initial and boundary conditions
\begin{align}
& U^n(\mathbf{x})=0,\quad \textrm{$\mathbf{x} \in \partial \Omega$, \ $n=2,3,\cdots, N$,}
\label{BCtime}
\\[3pt]
& U^0(\mathbf{x})=u_0(\mathbf{x}), \quad U^1(\mathbf{x})=u_0(\mathbf{x})+\tau u_1(\mathbf{x})+\frac{\tau^2}{2}\left(\Delta u_0-\mathbf{i}u_1-|u_0|^2u_0-wu_0\right),\quad \textrm{$\mathbf{x} \in \Omega$}.
\label{INITtime}
\end{align}
With the solutions of time discrete system \eqref{time-discrete}--\eqref{INITtime}, we split the errors into two parts
$$\| u^n-U^n_h\|_{L^2}\leq \| u^n-U^n\|_{L^2}+\| U^n-U^n_h\|_{L^2}.$$

Under this splitting, we will prove that the first term in the right-hand side of above inequality is bounded by $O(\tau^2)$ and the second term is bounded by $O(h^{2})$, with which, the classical inverse inequality and inductive assumption, we can obtain that the FE solutions are bounded in $L^{\infty}$--norm.   For the simplicity of notations, we use $C$ to denote a generic positive
constant and use  $\epsilon$ to denote	 a generic small positive constant, where $C,\epsilon$ are independent of spatial and temporal meshsize.

Next,  we give the following inequality, which will be used frequently.
\begin{lemma}
\label{gronwall}
{\it
Discrete Gronwall's inequality {\cite{Heywood,Nirenberg}} :
Let $\tau$, $B$ and $a_{k}$, $b_{k}$, $c_{k}$, $\gamma_{k}$,
for integers $k \geq 0$, be non--negative numbers such that
\[
a_j + \tau \sum_{k=0}^j b_{k}
\leq \tau \sum_{k=0}^j \gamma_{k} a_{k} +
\tau \sum_{k=0}^j c_{k} + B \, , \quad \mathrm{for } \quad j \geq 0 \, ,
\]
suppose that $\tau \gamma_{k} < 1$, for all $k$, and set $\sigma_{k}=(1-\tau \gamma_{k})^{-1}$. Then
\[
a_j + \tau \sum_{k=0}^j b_{k}
\leq  \exp(\tau \sum_{k=0}^j \gamma_{k} \sigma_{k}) (\tau \sum_{k=0}^j c_{k} + B) \, ,
\quad \mathrm{for } \quad j \geq 0 \, .
\]
}
\end{lemma}

\section{Temporal error estimates}\label{sec3}

In this subsection, we analyze the uniform boundedness of time--discrete solution $U^n$ in strong norms. Moreover, the error estimates of solution $U^n$ in certain norms are presented.

Let $u$ be the solution of the system (\ref{Schrodinger})--(\ref{INIT}), then $u^n$ satisfies
\begin{align}\label{truncation}
\delta^2_{\tau}u^n-\Delta u^{\overline{n}}+\mathbf{i}D_{\tau}u^n+|u^n|^2u^{\overline{n}}+wu^{\overline{n}}=P^n, \ n=1,2,\cdots, N-1,
\end{align}
where
$$P^n=(\delta^2_{\tau}u^n-u^n_{tt})-\Delta(u^{\overline{n}}-u^n)+\mathbf{i}(D_{\tau}u^n-u^n_t)+\left(|u^n|^2+w\right)\left(u^{\overline{n}}-u^n\right).$$
By using Taylor's expansion and the assumption (\ref{assumption}), one can obtain that
\begin{align}\label{estimate}
\left(\tau\sum^n_{k=1}\| P^k\|^2_{L^2}\right)^{\frac{1}{2}}+\left(\tau\sum^n_{k=1}\| \nabla P^k\|^2_{L^2}\right)^{\frac{1}{2}}\leq C\tau^2, \ n=1,2,\cdots, N-1.
\end{align}


\begin{theorem}\label{theorem2}
Suppose that the system (\ref{Schrodinger})--(\ref{INIT})  has a unique solution u satisfying (\ref{assumption}).
Then the time discrete system defined in (\ref{time-discrete})--(\ref{INITtime}) has a unique solution $U^m, m=1,\cdots,N$. Moreover, there exists $\tau_0>0$, such that when $\tau<\tau_0$,
\begin{align}
&\| U^m\|_{H^2}\leq M+1,\label{estimate1}\\
&\| e^m\|_{L^2}+\| \nabla e^m\|_{L^2}+\| \Delta e^m\|_{L^2}\leq C'\tau^2, \label{estimate2}
\end{align}
where $e^m=u^m-U^m$, $C'$ is a positive constant independent of $\tau$ and $h$.
\end{theorem}
\noindent{\bf Proof.} Since the system (\ref{time-discrete}) is a linear elliptic equation, we can get the existence and uniqueness of the solution $U^m$, $m=2,3,\cdots,N.$ First, we prove that there exists $\tau'_0>0$ such that the estimate (\ref{estimate2}) holds for $m=1,\cdots,N.$

Obviously, when $m=0$, $e^0=0.$  From \eqref{assumption} and \eqref{INITtime}, it is easy to see that, when $\tau\le \tau_1=\frac{C'}{C_1}$,
\begin{align}\label{e1-error}
\|e^1\|_{L^2}+\|\nabla e^1\|_{L^2}+\|\Delta e^1\|_{L^2} \le C_1\tau^3\le C'\tau^2.
\end{align}
Hence, \eqref{estimate2} holds for $m=1$.

Now suppose that \eqref{estimate2} is valid for $m\leq k-1$. Then there exists a positive constant $\tau_2$ such that for  $\tau\le\tau_2=\frac{1}{\sqrt{C_0C'}}$,
\begin{align}\label{en-H2-norm}
\|e^k\|_{L^{\infty}}\le C_{\infty}\|e^k\|_{H^2} \leq C_{\infty}C_0\left(\|e^k\|_{L^2}+\|\nabla e^k\|_{L^2}+\|\Delta e^k\|_{L^2}\right) \leq C_{\infty}C_0C'\tau^2\le C_{\infty}.
\end{align}

Next, we prove \eqref{estimate2} holds for $m=k$.
Subtracting Eq. \eqref{time-discrete} from Eq. \eqref{truncation} gives
\begin{align}\label{dalta-en-err-equ}
\delta^2_{\tau}e^n-\Delta e^{\overline{n}}+{\bf i}D_{\tau}e^n+R^n=P^n, \qquad 1\le n\le k-1,
\end{align}
where $R^n=(|u^n|^2+w)u^{\overline{n}}-(|U^n|^2+w)U^{\overline{n}}$.

Since
\begin{align}
R^n&=(|u^n|^2+w)u^{\overline{n}}-(|U^n|^2+w)U^{\overline{n}}\nonumber\\
&=|u^n|^2u^{\overline{n}}-|U^n|^2U^{\overline{n}}+we^{\overline{n}}\nonumber\\
&=u^n(u^n)^{*}u^{\bar{n}}-(u^n-e^n)((u^n)^{*}-(e^n)^*)(u^{\bar{n}}-e^{\bar{n}})+we^{\overline{n}}\nonumber\\
&=(u^n)^{*}u^{\bar{n}}e^n+u^nu^{\bar{n}}(e^n)^*+u^n(u^n)^{*}e^{\bar{n}}-(e^n)^*e^{\bar{n}}u^n-e^ne^{\bar{n}}(u^n)^{*}-e^n(e^n)^*u^{\bar{n}}+e^n(e^n)^*e^{\bar{n}}+we^{\overline{n}},\nonumber
\end{align}
one has
\begin{align}\label{Rn-L2-norm}
\|R^n\|_{L^2}=&\|(u^n)^{*}u^{\bar{n}}e^n+u^nu^{\bar{n}}(e^n)^*+u^n(u^n)^{*}e^{\bar{n}}-(e^n)^*e^{\bar{n}}u^n-e^ne^{\bar{n}}(u^n)^{*}-e^n(e^n)^*u^{\bar{n}}+e^n(e^n)^*e^{\bar{n}}+we^{\overline{n}}\|_{L^2}\nonumber\\
\leq&\|u^n\|_{L^{\infty}}\|u^{\bar{n}}\|_{L^{\infty}}\|e^n\|_{L^2}+\|u^n\|_{L^{\infty}}\|u^{\bar{n}}\|_{L^{\infty}}\|e^n\|_{L^2}+\|u^n\|^2_{L^{\infty}}\|e^{\bar{n}}\|_{L^2}
+\|e^n\|_{L^{\infty}}\|e^{\bar{n}}\|_{L^{2}}\|u^n\|_{L^\infty}\nonumber\\&+\|e^n\|_{L^{\infty}}\|e^{\bar{n}}\|_{L^{2}}\|u^n\|_{L^\infty}+\|e^n\|_{L^{\infty}}\|e^n\|_{L^2}\|u^{\bar{n}}\|_{L^\infty}+\|e^n\|^2_{L^{\infty}}\|e^{\bar{n}}\|_{L^2}+\|w\|_{L^{\infty}}\|e^{\bar{n}}\|_{L^2}\nonumber\\
\leq&C(\|e^{n-1}\|_{L^2}+\|e^n\|_{L^2}+\|e^{n+1}\|_{L^2})
\end{align}
where (\ref{assumption}), (\ref{en-H2-norm}) and the regularity of $w$ are used.

In addition,
\begin{align}
\nabla R^n&=\nabla((u^n)^*u^{\bar{n}}e^n)+\nabla (u^nu^{\bar{n}}(e^n)^*)+\nabla (u^n(u^n)^*e^{\bar{n}})-\nabla ((e^n)^*e^{\bar{n}}u^n)-\nabla (e^ne^{\bar{n}}(u^n)^*)\nonumber\\
&\;\;\;-\nabla(e^n(e^n)^*u^{\bar{n}})+\nabla(e^n(e^n)^*e^{\bar{n}})+\nabla(we^{\bar{n}})\nonumber\\
&\triangleq \sum^{8}_{j=1}I_j.
\end{align}
By using (\ref{assumption}) and (\ref{en-H2-norm}), one has
\begin{align}
\|I_1\|_{L^2}=&\|\nabla (u^n)^*u^{\bar{n}}e^n+(u^n)^*\nabla u^{\bar{n}}e^n+(u^n)^*u^{\bar{n}}\nabla e^n\|_{L^2}\nonumber\\
\leq& \|\nabla u^n\|_{L^{\infty}}\|u^{\bar{n}}\|_{L^{\infty}}\|e^n\|_{L^2}+\|u^n\|_{L^{\infty}}\|\nabla u^{\bar{n}}\|_{L^{\infty}}\|e^n\|_{L^2}+\|u^n||_{L^{\infty}}\|u^{\bar{n}}\|_{L^{\infty}}|\|\nabla e^n\|_{L^2}\nonumber\\
\leq& C(\|e^n\|_{L^2}+\|\nabla e^n\|_{L^2}),\label{2eq1}\\
\|I_2\|_{L^2}=&\|\nabla u^nu^{\bar{n}}(e^n)^*+u^n\nabla u^{\bar{n}}(e^n)^*+u^nu^{\bar{n}}\nabla (e^n)^*\|_{L^2}\nonumber\\
\leq& \|\nabla u^n\|_{L^{\infty}}\|u^{\bar{n}}\|_{L^{\infty}}\|e^n\|_{L^2}+\|u^n\|_{L^{\infty}}\|\nabla u^{\bar{n}}\|_{L^{\infty}}\|e^n\|_{L^2}+\|u^n\|_{L^{\infty}}\|u^{\bar{n}}\|_{L^{\infty}}\|\nabla e^n\|_{L^2}\nonumber\\
\leq& C(\|e^n\|_{L^2}+\|\nabla e^n\|_{L^2}),\label{2eq2}\\
\|I_3\|_{L^2}=&\|\nabla u^n(u^n)^*e^{\bar{n}}+ u^n\nabla(u^n)^*e^{\bar{n}}+u^n(u^n)^*\nabla e^{\bar{n}}\|_{L^2}\nonumber\\
\leq& 2\|\nabla u^n\|_{L^{\infty}}\|u^{n}\|_{L^{\infty}}\|e^{\bar{n}}\|_{L^2}+\|u^n\|^2_{L^{\infty}}\|\nabla e^{\bar{n}}\|_{L^2}\nonumber\\
\leq& C(\|e^{\bar{n}}\|_{L^2}+\|\nabla e^{\bar{n}}\|_{L^2})\nonumber\\
\leq& C(\|e^{n+1}\|_{L^2}+\|e^{n-1}\|_{L^2}+\|\nabla e^{n+1}\|_{L^2}+\|\nabla e^{n-1}\|_{L^2}),\label{2eq3}\\
\|I_4\|_{L^2}=&\|\nabla u^n(e^n)^*e^{\bar{n}}+u^n\nabla (e^n)^*e^{\bar{n}}+u^n(e^n)^*\nabla e^{\bar{n}}\|_{L^2}\nonumber\\
\leq& \|\nabla u^n\|_{L^{\infty}}\|e^{\bar{n}}\|_{L^2}\|e^n\|_{L^\infty}+\|u^n\|_{L^{\infty}}\|\nabla e^n\|_{L^2}\|e^{\bar{n}}\|_{L^\infty}+\|u^n\|_{L^{\infty}}\|e^n\|_{L^\infty}\|\nabla e^{\bar{n}}\|_{L^2}\nonumber\\
\leq& \|\nabla u^n\|_{L^{\infty}}\|e^n\|_{L^\infty}\|e^{\bar{n}}\|_{L^2}+C_{\infty}\|u^n\|_{L^{\infty}}\|e^n\|_{H^2}\| e^{\bar{n}}\|_{H^2}+\|u^n\|_{L^{\infty}}\|e^n\|_{L^{\infty}}\|\nabla e^{\bar{n}}\|_{L^2}\nonumber\\
\leq& \|\nabla u^n\|_{L^{\infty}}\|e^n\|_{L^\infty}\|e^{\bar{n}}\|_{L^2}+C_{\infty}C_0\|u^n\|_{L^{\infty}}\|e^n\|_{H^2}(\| e^{\bar{n}}\|+\| \nabla e^{\bar{n}}\|+\| \Delta e^{\bar{n}}\|)\nonumber\\
&+\|u^n\|_{L^{\infty}}\|e^n\|_{L^{\infty}}\|\nabla e^{\bar{n}}\|_{L^2}\nonumber\\
\leq& C(\|e^{\bar{n}}\|_{L^2}+\|\nabla e^{\bar{n}}\|_{L^2}+\|\Delta e^{\bar{n}}\|_{L^2})\nonumber\\
\leq&C(\|e^{n+1}\|_{L^2}+\|e^{n-1}\|_{L^2}+\|\nabla e^{n+1}\|_{L^2}+\|\nabla e^{n-1}\|_{L^2}+\|\Delta e^{n+1}\|_{L^2}+\|\Delta e^{n-1}\|_{L^2}),\label{2eq4}\\
\|I_5\|_{L^2}=&\|\nabla(u^n)^*e^{n}e^{\bar{n}}+(u^n)^*\nabla e^{n}e^{\bar{n}}+(u^n)^*e^{n}\nabla e^{\bar{n}}\|_{L^2}\nonumber\\
\leq& \|\nabla u^n\|_{L^{\infty}}\|e^{\bar{n}}\|_{L^2}\|e^n\|_{L^\infty}+\|u^n\|_{L^{\infty}}\|\nabla e^n\|_{L^2}\|e^{\bar{n}}\|_{L^\infty}+\|u^n\|_{L^{\infty}}\|e^n\|_{L^\infty}\|\nabla e^{\bar{n}}\|_{L^2}\nonumber\\
\leq& \|\nabla u^n\|_{L^{\infty}}\|e^{\bar{n}}\|_{L^2}\|e^n\|_{L^\infty}+C_{\infty}\|u^n\|_{L^{\infty}}\| e^n\|_{H^2}\|e^{\bar{n}}\|_{H^2}+\|u^n\|_{L^{\infty}}\|e^n\|_{L^\infty}\|\nabla e^{\bar{n}}\|_{L^2}\nonumber\\
\leq&C(\|e^{n+1}\|_{L^2}+\|e^{n-1}\|_{L^2}+\|\nabla e^{n+1}\|_{L^2}+\|\nabla e^{n-1}\|_{L^2}+\|\Delta e^{n+1}\|_{L^2}+\|\Delta e^{n-1}\|_{L^2}),\label{2eq5}
\end{align}
where the analysis here is the same as (\ref{2eq4}).
\begin{align}
\|I_6\|_{L^2}=&\|\nabla u^{\bar{n}}(e^n)^*e^n+u^{\bar{n}}\nabla (e^n)^*e^n+u^{\bar{n}}(e^n)^*\nabla e^n\|_{L^2}\nonumber\\
\leq& \|\nabla u^{\bar{n}}\|_{L^{\infty}}\|e^n\|_{L^\infty}\|e^n\|_{L^2}+2\|u^{\bar{n}}\|_{L^{\infty}}\|e^n\|_{L^\infty}\|\nabla e^n\|_{L^2}\nonumber\\
\leq& C(\|e^{n}\|_{L^2}+\|\nabla e^{n}\|_{L^2}),\label{2eq6}\\
\|I_7\|_{L^2}=&\|\nabla (e^n)^*e^ne^{\bar{n}}+(e^n)^*\nabla e^ne^{\bar{n}}+(e^n)^*e^n\nabla e^{\bar{n}}\|_{L^2}\nonumber\\
\leq &2\|e^n\|_{L^\infty}\|\nabla e^n\|_{L^2} \|e^{\bar{n}}\|_{L^\infty}+\|e^n\|^2_{L^\infty}\|\nabla e^{\bar{n}}\|_{L^2}\nonumber\\
\leq&2C_{\infty}\|e^n\|_{L^\infty}\|e^n\|_{H^2} \|e^{\bar{n}}\|_{H^2}+\|e^n\|^2_{L^\infty}\|\nabla e^{\bar{n}}\|_{L^2}\nonumber\\
\leq& C(\|e^{\bar{n}}\|_{L^2}+\|\nabla e^{\bar{n}}\|_{L^2}+\|\Delta e^{\bar{n}}\|_{L^2})\nonumber\\
\leq&C(\|e^{n+1}\|_{L^2}+\|e^{n-1}\|_{L^2}+\|\nabla e^{n+1}\|_{L^2}+\|\nabla e^{n-1}\|_{L^2}+\|\Delta e^{n+1}\|_{L^2}+\|\Delta e^{n-1}\|_{L^2}),\label{2eq7}\\
\|I_8\|_{L^2}=&\|\nabla we^{\bar{n}}+w\nabla e^{\bar{n}}\|_{L^2}\nonumber\\
\leq&\|\nabla w\|_{L^{\infty}}\|e^{\bar{n}}\|_{L^2}+\|w\|_{L^{\infty}}\|\nabla e^{\bar{n}}\|_{L^2}\nonumber\\
\leq& C(\|e^{n-1}\|_{L^2}+\|e^{n+1}\|_{L^2}+\|\nabla e^{n-1}\|_{L^2}+\|\nabla e^{n+1}\|_{L^2}).\label{2eq8}
\end{align}

Adding the equations from (\ref{2eq1}) to (\ref{2eq8}), one has
\begin{multline}\label{Rn-H1-norm}
\|\nabla R^n\|_{L^2}\le C(\|e^{n-1}\|_{L^2}+\|e^{n}\|_{L^2}+\|e^{n+1}\|_{L^2}+\|\nabla e^{n-1}\|_{L^2}+\|\nabla e^{n}\|_{L^2}+\|\nabla e^{n+1}\|_{L^2}+\|\Delta e^{n-1}\|_{L^2}+\|\Delta e^{n+1}\|_{L^2}),
\end{multline}
for $1\le n\le k-1$.

Computing the inner product with $D_{\tau} e^{n}$ on both sides of \eqref{dalta-en-err-equ} and taking the real parts, one gets
\begin{align}\label{error18}
&\frac{1}{2\tau}\left(\left\|\frac{e^{n+1}-e^n}{\tau}\right\|^2_{L^2}-\left\|\frac{e^{n}-e^{n-1}}{\tau}\right\|^2_{L^2}\right)+\frac{\|\nabla e^{n+1}\|^2_{L^2}-\|\nabla e^{n-1}\|^2_{L^2}}{4\tau}\nonumber\\
=&{\rm Re}\left(P^n-R^n,\frac{e^{n+1}-e^{n-1}}{2\tau}\right)\nonumber\\
\leq& \frac{1}{4}\left(\left\|\frac{e^{n+1}-e^{n}}{\tau}\right\|^2_{L^2}+\left\|\frac{e^{n}-e^{n-1}}{\tau}\right\|^2_{L^2}\right)+\|  P^n\|^2_{L^2}+\| R^n\|^2_{L^2}.
\end{align}
Multiplying above inequality both sides by $2\tau$ and summing up from $n=1$ to $n=k-1$, one obtains
\begin{align}\label{error19}
&\left\|\frac{e^{k}-e^{k-1}}{\tau}\right\|^2_{L^2}+\frac{\|\nabla e^{k}\|^2_{L^2}+\|\nabla e^{k-1}\|^2_{L^2}}{2}\nonumber\\
\le&\frac{\|\nabla e^1\|^2_{L^2}+\|\nabla e^0\|^2_{L^2}}{2}+\left\|\frac{e^1-e^0}{\tau}\right\|^2_{L^2}+\tau\sum^{k-1}_{n=0}\left\|\frac{e^{n+1}-e^n}{\tau}\right\|^2_{L^2}+2\tau\sum^{k-1}_{n=1}(\|P^n\|^2_{L^2}+\|R^n\|^2_{L^2})\nonumber\\
\le& \tau\sum^{k-1}_{n=0}\left\|\frac{e^{n+1}-e^{n}}{\tau}\right\|^2_{L^2}+C\tau\sum^{k-1}_{n=0}\|e^{n+1}\|^2_{L^2}+C\tau^4.\qquad\qquad\,\mbox{(use \eqref{estimate}, \eqref{e1-error} and \eqref{Rn-L2-norm})}
\end{align}

%

Now computing the  inner product with $-\Delta D_{\tau} e^{n}$ on both sides of (\ref{dalta-en-err-equ}) and taking the real parts, one has
\begin{align}\label{error24}
&\frac{1}{2\tau}\left(\left\|\frac{\nabla(e^{n+1}-e^n)}{\tau}\right\|^2_{L^2}-\left\|\frac{\nabla(e^n-e^{n-1})}{\tau}\right\|^2_{L^2}\right)+\frac{\|\Delta e^{n+1}\|^2_{L^2}-\|\Delta e^{n-1}\|^2_{L^2}}{4\tau}\nonumber\\
=&{\rm Re}\left(P^n-R^n,\frac{-\Delta(e^{n+1}-e^{n-1})}{2\tau}\right)\nonumber\\
\leq& \frac{1}{4}\left(\left\|\frac{\nabla(e^{n+1}-e^n)}{\tau}\right\|^2_{L^2}+\left\|\frac{\nabla(e^n-e^{n-1})}{\tau}\right\|^2_{L^2}\right)+\| \nabla P^n\|^2_{L^2}+\| \nabla R^n\|^2_{L^2}.
\end{align}
Multiplying above inequality both sides by $2\tau$ and summing up from $n=1$ to $n=k-1$, one gets,
\begin{align}\label{error25}
&\left\|\frac{\nabla(e^{k}-e^{k-1})}{\tau}\right\|^2_{L^2}+\frac{\|\Delta e^{k}\|^2_{L^2}+\|\Delta e^{k-1}\|^2_{L^2}}{2}\nonumber\\
\leq&
\frac{\|\Delta e^1\|^2_{L^2}+\|\Delta e^0\|^2_{L^2}}{2}+\left\|\frac{\nabla(e^1-e^0)}{\tau}\right\|^2_{L^2}+\tau\sum^{k-1}_{n=0}\left\|\frac{\nabla(e^{n+1}-e^n)}{\tau}\right\|^2_{L^2}\nonumber\\
&+2\tau\sum^{k-1}_{n=1}(\|\nabla P^n\|^2_{L^2}+\|\nabla R^n\|^2_{L^2})\nonumber\\
\le&
\tau\sum^{k-1}_{n=0}\left\|\frac{\nabla(e^{n+1}-e^{n})}{\tau}\right\|^2_{L^2}+C\tau\sum^{k-1}_{n=0}(\|e^{n+1}\|^2_{L^2}+\|\nabla e^{n+1}\|^2_{L^2}\nonumber\\
&+\|\Delta e^{n+1}\|^2_{L^2})+C\tau^4. \qquad\mbox{(use \eqref{estimate}, \eqref{e1-error} and \eqref{Rn-H1-norm})}
\end{align}

Adding Eqs. \eqref{error19} and \eqref{error25}, we get
\begin{align}\label{error00}
&\left\|\frac{e^{k}-e^{k-1}}{\tau}\right\|^2_{L^2}+\frac{\|\nabla e^{k}\|^2_{L^2}+\|\nabla e^{k-1}\|^2_{L^2}}{2}+ \left\|\frac{\nabla(e^{k}-e^{k-1})}{\tau}\right\|^2_{L^2}+\frac{\|\Delta e^{k}\|^2_{L^2}+\|\Delta e^{k-1}\|^2_{L^2}}{2}\nonumber\\
\le& \tau\sum^{k-1}_{n=0}\left(\left\|\frac{\nabla(e^{n+1}-e^{n})}{\tau}\right\|^2_{L^2}+\left\|\frac{e^{n+1}-e^{n}}{\tau}\right\|^2_{L^2}\right)\nonumber\\
&+C\tau\sum^{k-1}_{n=0}(\|e^{n+1}\|^2_{L^2}+\|\nabla e^{n+1}\|^2_{L^2}+\|\Delta e^{n+1}\|^2_{L^2})+C\tau^4 \nonumber \\
\le& C\tau\sum^{k}_{n=1}\bigg(\left\|\frac{e^{n}-e^{n-1}}{\tau}\right\|^2_{L^2}+\frac{\|\nabla e^{n}\|^2_{L^2}+\|\nabla e^{n-1}\|^2_{L^2}}{2}\nonumber\\
&+\left\|\frac{\nabla(e^{n}-e^{n-1})}{\tau}\right\|^2_{L^2}+\frac{\|\Delta e^{n}\|^2_{L^2}+\|\Delta e^{n-1}\|^2_{L^2}}{2}\bigg)+C\tau^4,
\end{align}
where we have used the Poincar$\acute{\rm{e}}$ inequality $\|e^{n}\|_{L^2} \le C \|\nabla e^{n}\|_{L^2}$.

By using Gronwall's inequality to \eqref{error00}, one obtains that, there exists $\tau_3>0$, such that for $\tau\le\tau_3$,
\begin{align}\label{error26}
\|\nabla e^k\|^2_{L^2} + \|\Delta e^{k}\|^2_{L^2}\leq C \tau^4,
\end{align}
which implies that
\begin{align}\label{error27}
\| e^k\|_{L^2} + \|\nabla e^k\|_{L^2} + \|\Delta e^{k}\|_{L^2}\leq C_2 \tau^2.
\end{align}
Now we take $C'=C_2$. After $C'$ is obtained, we let $\tau'_0=\min\{\tau_1, \tau_2,\tau_3\}$, then \eqref{estimate2} holds for $m=k$. The mathematical induction of \eqref{estimate2} is finished. Therefore, \eqref{estimate2} holds for $m=1,2,...,N.$
Thus, there exists a positive constant $\tau_4$ such that for $\tau\le\tau_4=\frac{1}{\sqrt{C'C_0}}$,
\begin{align*}
&\|e^m\|_{H^2}\le C_0(\|e^m\|_{L^2}  + \|\nabla e^m\|_{L^2} + \|\Delta e^{m}\|_{L^2}) \le C'C_0\tau^2 \le 1,\\
&\|U^m\|_{H^2}\le \|u^m\|_{H^2}+\|e^m\|_{H^2}\le M+1,
\end{align*}
for $m=1,2,...,N$. Let $\tau_0=\min\{\tau'_0, \tau_4\}$, we finish the proof of \eqref{estimate1} and \eqref{estimate2}.\qed

\section{Spatial error estimates}\label{sec4}

In this section, we will study the uniform boundedness of solution $U^n_h$ in strong norms.

\begin{lemma}
{\it
Suppose that the time--discrete system \emph{(\ref{time-discrete})--\eqref{INITtime}} has a unique solution $U^n$ \emph{(}$n=0,\cdots,N$\emph{)}, then
\begin{align}\label{bound}
\| R_hU^n\|_{L^{\infty}}\leq M_1, \quad  n=0,1,\cdots,N.
\end{align}
}
\end{lemma}
{\bf Proof.}
By using \eqref{estimate1}  and the Sobolev embedding theorem, we have
\begin{align}\label{Rhinfty-norm}
\max\limits_{0\le n\le N}\| R_hU^n\|_{L^{\infty}}\leq&C\max\limits_{0\le n\le N}\| R_hU^n\|_{W^{1,6}}\nonumber\\
\leq&C\max\limits_{0\le n\le N}\|U^n\|_{W^{1,6}}\nonumber \qquad\mbox{(use \eqref{theory2})} \\
\leq&C\max\limits_{0\le n\le N}\|U^n\|_{H^2}\nonumber\\
\leq& M_1.
\end{align}
Thus, the proof of this lemma is complete. \qed

The variation form of (\ref{time-discrete}) can be defined by
\begin{align}\label{error31}
(\delta^2_{\tau}U^n, v)+(\nabla U^{\overline{n}},\nabla v)+\mathbf{i}(D_{\tau}U^n,v)+\left(|U^n|^2U^{\overline{n}},v\right)+(wU^{\overline{n}},v)=0,\quad \forall v\in H^1_0(\Omega),
\end{align}
for $n=1,2,\cdots,N-1$.

Subtracting (\ref{FEMscheme}) from (\ref{error31}) and using (\ref{Ritz}), one gets
\begin{multline}\label{error32}
(\delta^2_{\tau}e^n_h, v_h)+(\nabla e_h^{\overline{n}},\nabla v_h)+\mathbf{i}(D_{\tau}e^n_h,v_h)+\left(Q^n,v_h\right)=-(\delta^2_{\tau}(U^n-R_hU^n), v_h)\\+\mathbf{i}(D_{\tau}(U^n-R_hU^n),v_h), \;\;n=1, 2,\cdots,N-1,\quad \forall v_h\in V_h,
\end{multline}
where $Q^n=\left(|U^n|^2+w\right)U^{\overline{n}}-\left(|U^n_h|^2+w\right)U^{\overline{n}}_h$ and $e^n_h=R_hU^n-U^n_h$.

\begin{theorem}\label{theorem3}
Suppose that the system \emph{(\ref{Schrodinger})--(\ref{INIT})}  has a unique solution $u$ satisfying \emph{(\ref{assumption})}.
Then the finite element equation defined by \emph{(\ref{FEMscheme})} has a unique solution $U^n_h, n=2,\cdots,N$.  Moreover, there exists $\tau^*>0, h^*>0$, such that,
\begin{align}
\left\| U^n_h\right\|_{L^{\infty}}\leq M_1+1, \label{estimate3}
\end{align}
 when $\tau<\tau^*, h<h^*$.
\end{theorem}
{\bf Proof.}  Firstly, we prove the existence and uniqueness of solution $U^n_h$ of system (\ref{FEMscheme}). If the solution $U^n_h$ is given for $n=0,1,2,\cdots,m-1$, then the system (\ref{FEMscheme}) has a unique solution if and only if the following homogeneous equation
\begin{align}\label{Wh-homogeneous-equ}
\frac{1}{\tau^2}(W_h, v_h)+\frac{1}{2}(\nabla W_h, \nabla v_h)+\frac{\bf i}{2\tau}(W_h, v_h)+\frac{1}{2}(|U^{m-1}_h|^2W_h, v_h)+\frac{1}{2}(w({\bf x})W_h, v_h)=0, \quad \forall v_h\in V_h,
\end{align}
has only zero solution.

Let $v_h=W_h$ in above equation and take the imaginary part, one has
\begin{align*}
\|W_h\|_{L^2}=0,
\end{align*}
which implies that $W_h=0$. Thus, equation \eqref{Wh-homogeneous-equ} has only zero solution. It is natural to obtain that the uniqueness and existence for the solutions of (\ref{FEMscheme}).

Before the proof of \eqref{estimate3} is given, we study the following result: there exists two constants $\tau^*>0$ and $h^*>0$, such that, when $\tau\le\tau^*$ and $h\le h^*$,
\begin{align}\label{primary-result-1}
\|e^k_h\|_{L^2}\le h^{\frac{5}{3}},\quad 0\le k\le N.
\end{align}
As $k=0, 1$, from \eqref{INITtime} and \eqref{FEMIN}, we can get $e^0_h=0$, $e^1_h=0$. Thus, \eqref{primary-result-1} holds for $k=0, 1$, and
\begin{align}\label{Unh-Linfty-0,1}
\max\limits_{0\le k\le 1}\|U^k_h\|_{L^{\infty}}\le& \max\limits_{0\le k\le 1}\|R_hU^k\|_{L^{\infty}}+\max\limits_{0\le k\le 1}\|e^k_h\|_{L^{\infty}}\nonumber\\
\le& M_1+C_3h^{-\frac{d}{2}}\max\limits_{0\le k\le 1}\|e^k_h\|_{L^2}\qquad\mbox{({\color{blue}use \eqref{Rhinfty-norm} and {\eqref{inversein}}})}       \nonumber\\
\le& M_1+C_3h^{-\frac{d}{2}}h^{\frac{5}{3}}\qquad\qquad\quad\,\,\,\mbox{(use \eqref{primary-result-1})}       \nonumber\\
\le& M_1+1,
\end{align}
\vspace{-0.2in}\\
when $h\le h_1=\left(\frac{1}{C_3}\right)^{\frac{6}{10-3d}}$ ($d=2,3$)\,($h_1>0$), where the inverse inequality $\|e_h\|_{L^{\infty}}\le C_3h^{-\frac{d}{2}}|e_h\|_{L^{2}}$ is used.

Now suppose that (\ref{primary-result-1}) is valid for $k\leq m-1$,
we can get
\begin{align}\label{Ukh-infty-assum}
\max_{0\le n\le m-1}\|U^n_h\|_{L^{\infty}}\le& \max_{0\le n\le m-1}\|e^n_h\|_{L^{\infty}}+\max_{0\le n\le m-1}\|R_hU^n\|_{L^{\infty}}\nonumber\\
\le&C_3h^{-\frac{d}{2}}\max_{0\le n\le m-1}\|e^n_h\|_{L^2}+M_1\nonumber\\
\le& 1+M_1,
\end{align}
when $h\le h_1$. Next, we want to show that (\ref{primary-result-1}) is also valid for $k=m$.

Putting $v_h=D_{\tau}e^n_h$($1\le n \le m-1$) on both sides of (\ref{error32})  and taking the real parts, one has
\begin{align}\label{error39}
&\frac{1}{2\tau}\left(\left\|\frac{e^{n+1}_h-e^n_h}{\tau}\right\|^2_{L^2}-\left\|\frac{e^{n}_h-e^{n-1}_h}{\tau}\right\|^2_{L^2}\right)+\frac{1}{4\tau}\left(\left\|\nabla e^{n+1}_h\right\|^2_{L^2}-\left\|\nabla e^{n-1}_h\right\|^2_{L^2}\right)\nonumber\\
=&-Re\left(Q^n,D_{\tau}e^n_h\right)-Re(\delta^2_{\tau}(U^n-R_hU^n), D_{\tau}e^n_h)+Im(D_{\tau}(U^n-R_hU^n),D_{\tau}e^n_h)\nonumber\\
\leq&\| Q^n\|^2_{L^2}+\|\delta^2_{\tau}(U^n-R_hU^n)\|^2_{L^2}+\| D_{\tau}(U^n-R_hU^n)\|^2_{L^2}+ \left\|\frac{e^{n+1}_h-e^n_h}{\tau}\right\|^2_{L^2}+\left\|\frac{e^n_h-e^{n-1}_h}{\tau}\right\|^2_{L^2}\nonumber\\
\leq&\| Q^n\|^2_{L^2}+Ch^4\|\delta^2_{\tau}U^n\|^2_{H^2}+Ch^4\| D_{\tau}U^n\|^2_{H^2}+ \left\|\frac{e^{n+1}_h-e^n_h}{\tau}\right\|^2_{L^2}+\left\|\frac{e^n_h-e^{n-1}_h}{\tau}\right\|^2_{L^2}.
\end{align}
And $\| Q^n\|_{L^2}$ can be further analyzed, thus
\begin{align}\label{error40}
&\frac{1}{2\tau}\left(\left\|\frac{e^{n+1}_h-e^n_h}{\tau}\right\|^2_{L^2}-\left\|\frac{e^{n}_h-e^{n-1}_h}{\tau}\right\|^2_{L^2}\right)+\frac{1}{4\tau}\left(\left\|\nabla e^{n+1}_h\right\|^2_{L^2}-\left\|\nabla e^{n-1}_h\right\|^2_{L^2}\right)\nonumber\\
\le&C\|(|U^n|+|U^n_h|)(U^n-R_hU^n)U^{\overline{n}}\|^2_{L^2}+C\|(|U^n|+|U^n_h|)e^n_hU^{\overline{n}}\|^2_{L^2}+C\|(|U^n_h|^2+w)(U^{\overline{n}}-R_hU^{\overline{n}})\|^2_{L^2}\nonumber\\
&+C\|(|U^n_h|^2+w)e^{\overline{n}}_h\|^2_{L^2}+Ch^4(\|\delta^2_{\tau}e^n\|^2_{H^2}+\|\delta^2_{\tau}u^n\|^2_{H^2})+Ch^4(\| D_{\tau}e^n\|^2_{H^2}+\| D_{\tau}u^n\|^2_{H^2})\nonumber\\
&+C\left\|\frac{e^{n+1}_h-e^n_h}{\tau}\right\|^2_{L^2}+C\left\|\frac{e^n_h-e^{n-1}_h}{\tau}\right\|^2_{L^2}\nonumber\\
\le&
Ch^4\|U^n\|^2_{H^2}(\|U^n\|^2_{L^{\infty}}+\|U^n_h\|^2_{L^{\infty}})\|U^{\overline{n}}\|^2_{L^{\infty}}+C\|U^{\overline{n}}\|^2_{L^{\infty}}(\|U_h^n\|^2_{L^{\infty}}+\|U^{n}\|^2_{L^{\infty}})\|e^n_h\|^2_{L^2}\quad\,\,\,\mbox{(use \eqref{theory})}\nonumber\\
&+Ch^4\|U^{\overline{n}}\|^2_{H^2}(\|U^n_h\|^4_{L^{\infty}}+\|w\|^2_{L^{\infty}})+C(\|U^n_h\|^4_{L^{\infty}}+\|w\|^2_{L^{\infty}})\|e^{\overline{n}}_h\|^2_{L^2}  \qquad\qquad\qquad\mbox{(use \eqref{theory})}\nonumber\\
&+Ch^4+C\left\|\frac{e^{n+1}_h-e^n_h}{\tau}\right\|^2_{L^2}+C\left\|\frac{e^n_h-e^{n-1}_h}{\tau}\right\|^2_{L^2}\qquad\qquad\qquad\qquad\,\mbox{(use \eqref{assumption} and \eqref{estimate2})}\nonumber\\
\le&
C(\|e^n_h\|^2_{L^2}+\|e_h^{\overline{n}}\|^2_{L^2}+h^4)+C\left\|\frac{e^{n+1}_h-e^n_h}{\tau}\right\|^2_{L^2}+C\left\|\frac{e^n_h-e^{n-1}_h}{\tau}\right\|^2_{L^2}.\quad\,\,\mbox{(use \eqref{estimate1} and \eqref{Unh-Linfty-0,1})}
\end{align}

%

Summing above inequality from $n=1$ to $n=m-1$, we obtain
\begin{align*}
&\left\|\frac{e^{m}_h-e^{m-1}_h}{\tau}\right\|^2_{L^2}+\frac{\|\nabla e^{m}_h\|^2_{L^2}+\|\nabla e^{m-1}_h\|^2_{L^2}}{2}\nonumber\\
\le&
\left\|\frac{e^1_h-e^0_h}{\tau}\right\|^2_{L^2}+\frac{\|\nabla e^1_h\|^2_{L^2}+\|\nabla e^0_h\|^2_{L^2}}{2}+C\tau\sum^{m-1}_{n=0}\left\|\frac{e^{n+1}_h-e^n_h}{\tau}\right\|^2_{L^2}+C\tau\sum_{n=0}^{m}\|e^n_h\|^2_{L^2}+Ch^4\nonumber\\
\le&
C\tau\sum^{m-1}_{n=0}\left\|\frac{e^{n+1}_h-e^n_h}{\tau}\right\|^2_{L^2}+C\tau\sum^{m-1}_{n=0}\frac{\|\nabla e^{n+1}_h\|^2_{L^2}+\|\nabla e^n_h\|^2_{L^2}}{2}+Ch^4,
\end{align*}
where $e^0_h=0$, $e^1_h=0$ and the Poincar$\acute{\rm{e}}$ inequality $\|e^{n}_h\|_{L^2} \le C \|\nabla e^{n}_h\|_{L^2}$ are used.

Using the Gronwall's inequality, there is a $\tau_{5}>0$, when $\tau\le\tau_{5}$
\begin{align*}
\|\nabla e^m_h\|_{L^2}\le C_4h^2.
\end{align*}

By the fact that  $e^n_h\in H^1_0$, we have
\begin{align}
\|e^m_h\|_{L^2}\le C\|\nabla e^m_h\|_{L^2}\le C_5h^2.
\end{align}
There exists a positive constant $h_2=\left(\frac{1}{C_5}\right)^3>0$, such that, when $h\le h_2$,
\begin{align}
\|e^m_h\|_{L^2}\le C_5h^2\le h^{\frac{5}{3}}.
\end{align}
Thus, \eqref{primary-result-1} holds for $k=m$. Let $\tau^*=\tau_5$, $h^*=\min \{ h_1, h_2\}$, the induction proof for \eqref{primary-result-1} is closed.

Next, we prove the \eqref{estimate3} holds for $n=2,\dots,N$. From \eqref{primary-result-1}, we can get
\begin{align}
\max_{2\le n\le N}\|U^n_h\|_{L^{\infty}}
\le&
 \max_{2\le n\le N}\|R_hU^n\|_{L^{\infty}}+\max_{2\le n\le N}\|e^n_h\|_{L^{\infty}}\nonumber\\
\le&
\max_{2\le n\le N}\|R_hU^n\|_{L^{\infty}}+C_3h^{-\frac{d}{2}}\max_{2\le n\le N}\|e^n_h\|_{L^2}\nonumber\\
 \le&
\max_{2\le n\le N}\|R_hU^n\|_{L^{\infty}}+C_3h^{-\frac{d}{2}}h^{\frac{5}{3}}\nonumber\\
\le& M_1+1,
\end{align}
when $h\le h^*$.

The proof is complete.\qed

%


\section{Convergent analysis for the full discrete scheme}\label{sec5}

In this section, we  will provide the convergent analysis for the full discrete scheme (\ref{FEMscheme})--(\ref{FEMIN}).

\begin{theorem}\label{theorem1}
Suppose that the system \emph{(\ref{Schrodinger})--(\ref{INIT})}  has a unique solution u satisfying \emph{(\ref{assumption})}.
Then the finite element system defined in \emph{(\ref{FEMscheme})}--\emph{(\ref{FEMIN})} has a unique solution $U^n_h, n=0,\cdots,N$ and there exists $\tau_*>0$, $h_*>0$ such that
\begin{align}\label{fulles}
\| u^n-U^n_h\|_{L^2}\leq C^0(\tau^2+h^{r+1}),
\end{align}
when $\tau<\tau_*, h<h_*$, where $C^0$ is a positive constant independent of $\tau$ and $h$.
\end{theorem}

{\bf Proof.} As $n=0, 1$, from \eqref{assumption} and \eqref{theory}, we have
\begin{align}
\|u_0-U^0_h\|_{L^2}\le& \|u_0-R_hu_0\|_{L^2}\le Ch^{r+1}\|u_0\|_{H^{r+1}}\le Ch^{r+1},\nonumber\\
\|u^1-U^1_h\|_{L^2}=&\left\|u^1-R_h\left(u_0+\tau u_1+\frac{\tau^2}{2!}(\Delta u_0-{\bf i}u_1-|u_0|^2u_0-wu_0)\right)\right\|_{L^2}\nonumber\\
\le&
\bigg\|u_0+\tau u_1+\frac{\tau^2}{2!}(\Delta u_0-{\bf i}u_1-|u_0|^2u_0-wu_0)+\frac{\tau^3}{3!}(u_0)_{ttt}+O(\tau^3)\nonumber\\
&-R_h\left(u_0+\tau u_1+\frac{\tau^2}{2!}(\Delta u_0-{\bf i}u_1-|u_0|^2u_0-wu_0)\right)\bigg\|_{L^2}\nonumber\\
\le&
Ch^{r+1}\left\|u_0+\tau u_1+\frac{\tau^2}{2!}(\Delta u_0-{\bf i}u_1-|u_0|^2u_0-wu_0)\right\|_{H^{r+1}}+C\tau^3\nonumber\\
\le& Ch^{r+1}+C\tau^2.
\end{align}
Thus, \eqref{fulles} holds for $n=0,1$.

Next, we prove \eqref{fulles} holds for $2\le n\le N$. The weak form of system \eqref{truncation} can be written as:
\begin{align}\label{un-weak-form}
(\delta^2_{\tau}u^n, v)+(\nabla u^{\overline{n}}, \nabla v)+{\bf i}(D_{\tau}u^n, v)+(|u^n|^2u^{\overline{n}} ,v)+(wu^{\overline{n}}, v)=(P^n, v), \quad \forall v\in H^1_0,
\end{align}
where $1\le n\le N-1$.

Subtracting \eqref{un-weak-form} from \eqref{FEMscheme} and applying \eqref{Ritz}, we can get
\begin{align}\label{sigma-equ}
&(\delta^2_{\tau}\sigma^n_h, v_h)+(\nabla\sigma^{\overline{n}}_h, \nabla v_h)+{\bf i}(D_{\tau}\sigma^n_h, v_h)+(G^n, v_h)+(w\sigma^{\overline{n}}_h, v_h)\nonumber\\
=&(P^n, v_h)-(\delta^2_{\tau}(u^n-R_hu^n), v_h)-{\bf i}(D_{\tau}(u^n-R_hu^n), v_h)-(w(u^{\overline{n}}-R_hu^{\overline{n}}), v_h),
\end{align}
where $G^n=|u^n|^2u^{\overline{n}}-|U^n_h|^2U^{\overline{n}}_h$,\quad$\sigma^n_h=R_hu^n-U^n_h$.

From \eqref{theory}, \eqref{assumption} and \eqref{estimate3}, we can get
\begin{align}\label{Gnerr-estimate}
\|G^n\|_{L^2}=&\||u^n|^2u^{\overline{n}}-|U^n_h|^2U_h^{\overline{n}}\|_{L^2}\nonumber\\
\le&\||u^n|^2(u^{\overline{n}}-U_h^{\overline{n}})\|_{L^2}+\|(|u^n|^2-|U^n_h|^2)U_h^{\overline{n}}\|_{L^2}\nonumber\\
\le&\||u^n|^2(u^{\overline{n}}-R_hu^{\overline{n}})\|_{L^2}+\||u^n|^2\sigma^{\overline{n}}_h\|_{L^2}\nonumber\\
&+\|(u^n-R_hu^n)(|u^n|+|U^n_h|)U_h^{\overline{n}}\|_{L^2}+\|\sigma^n_h(|u^n|+|U^n_h|)U_h^{\overline{n}}\|_{L^2}\nonumber\\
\le& Ch^{r+1}\|u^{\overline{n}}\|_{H^{r+1}}\|u^n\|^2_{L^{\infty}}+C\|u^n\|^2_{L^{\infty}}(\|\sigma^{n-1}_h\|_{L^2}+\|\sigma^{n+1}_h\|_{L^2})\nonumber\\
&+Ch^{r+1}\|u^n\|_{H^{r+1}}(\|u^n\|_{L^\infty}+\|U^n_h\|_{L^\infty})\|U^{\overline{n}}_h\|_{L^{\infty}}\nonumber\\
&+C\|U^{\overline{n}}_h\|_{L^{\infty}}(\|u^n\|_{L^\infty}+\|U^n_h\|_{L^\infty})\|\sigma^n_h\|_{L^2}\nonumber\\
\le&
C(\|\sigma_h^{n-1}\|_{L^2}+\|\sigma_h^{n}\|_{L^2}+\|\sigma_h^{n+1}\|_{L^2})+Ch^{r+1}.
\end{align}

Putting $v_h=D_{\tau}\sigma^n_h$ in \eqref{sigma-equ} and taking the real parts, we get
\begin{align}
&\frac{1}{2\tau}\left(\left\|\frac{\sigma^{n+1}_h-\sigma^{n}_h}{\tau} \right\|^2_{L^2}-\left\|\frac{\sigma^{n}_h-\sigma^{n-1}_h}{\tau} \right\|^2_{L^2}\right)+\frac{\|\nabla\sigma^{n+1}_h\|^2_{L^2}-\|\nabla\sigma^{n-1}_h\|^2_{L^2}}{4\tau}\nonumber\\
\le&|(G^n, D_{\tau}\sigma^n_h)|+|(w\sigma^{\bar{n}}_h, D_{\tau}\sigma^n_h)|+|(\delta^2_{\tau}(u^n-R_hu^n), D_{\tau}\sigma^n_h)|+|(P^n, D_{\tau}\sigma^n_h)|\nonumber\\
&+|(D_{\tau}(u^n-R_hu^n), D_{\tau}\sigma^n_h)|+|(w(u^{\bar{n}}-R_hu^{\bar{n}}),D_{\tau}\sigma^n_h)|\nonumber\\
\le&C\|G^n\|^2_{L^2}+C\|w\|^2_{L^\infty}\|\sigma^{\bar{n}}_h\|^2_{L^2}+Ch^{2(r+1)}\|\delta^2_{\tau}u^n\|^2_{H^{r+1}}+C\|P^n\|^2_{L^2}\quad\,\mbox{(use \eqref{theory})}\nonumber\\
&+Ch^{2(r+1)}\|D_{\tau}u^n\|^2_{H^{r+1}}+Ch^{2(r+1)}\|w\|^2_{L^{\infty}}\|u^{\bar{n}}\|^2_{H^2}\qquad\qquad\quad\,\,\,\,\,\,\,\mbox{(use \eqref{theory})}\nonumber\\
&+C\left\|\frac{\sigma_h^{n+1}-\sigma_h^{n}}{\tau}\right\|_{L^2}^2+C\left\|\frac{\sigma_h^{n}-\sigma_h^{n-1}}{\tau}\right\|_{L^2}^2\nonumber\\
\le&C(\|\sigma^{n-1}_h\|^2_{L^2}+\|\sigma^{n}_h\|^2_{L^2}+\|\sigma^{n+1}_h\|^2_{L^2})+C\|P^n\|^2_{L^2}+Ch^{2(r+1)}\qquad\,{\mbox{(use \eqref{assumption} and \eqref{Gnerr-estimate})}}\nonumber\\
&+C\left\|\frac{\sigma_h^{n+1}-\sigma_h^{n}}{\tau}\right\|_{L^2}^2+C\left\|\frac{\sigma_h^{n}-\sigma_h^{n-1}}{\tau}\right\|_{L^2}^2.
\end{align}

Summing above inequality from $n=1$ to $n=N-1$, we have
\begin{align}\label{sigmaresult}
&\max_{1\le n\le N-1}\left\|\frac{\sigma_h^{n+1}-\sigma_h^{n}}{\tau}\right\|^2_{L^2}+\max_{1\le n\le N-1}\frac{\|\nabla\sigma^{n+1}_h\|^2_{L^2}+\|\nabla\sigma^{n}_h\|^2_{L^2}}{2}\nonumber\\
\le&\left\|\frac{\sigma_h^{1}-\sigma_h^{0}}{\tau}\right\|^2_{L^2}+\frac{\|\nabla\sigma^1_h\|^2_{L^2}+\|\nabla\sigma^0_h\|^2_{L^2}}{2}+C\tau\sum^N_{n=0}\|\sigma_h^n\|^2_{L^2}+C\tau\sum^{N-1}_{n=1}\|P^n\|^2_{L^2}\nonumber\\
&+C\tau\sum^{N-1}_{n=0}\left\|\frac{\sigma^{n+1}_h-\sigma^{n}_h}{\tau}\right\|^2_{L^2}+Ch^{2(r+1)}\nonumber\\
\le&C\tau\sum^N_{n=0}\|\sigma_h^n\|^2_{L^2}+C\tau\sum^{N-1}_{n=1}\|P^n\|^2_{L^2}+C\tau\sum^{N-1}_{n=0}\left\|\frac{\sigma^{n+1}_h-\sigma^{n}_h}{\tau}\right\|^2_{L^2}+C\tau^4+Ch^{2(r+1)}\nonumber\\
\le&C\tau\sum^N_{n=0}\|\sigma_h^n\|^2_{L^2}+C\tau\sum^{N-1}_{n=0}\left\|\frac{\sigma^{n+1}_h-\sigma^{n}_h}{\tau}\right\|^2_{L^2}+C\tau^4+Ch^{2(r+1)}\qquad\qquad\qquad\,\,\,\,{\mbox{(use \eqref{estimate})}}\nonumber\\
\le&C\tau\sum^{N-1}_{n=0}\frac{\|\nabla\sigma^{n+1}_h\|^2_{L^2}+\|\nabla\sigma^{n}_h\|^2_{L^2}}{2}+C\tau\sum^{N-1}_{n=0}\left\|\frac{\sigma^{n+1}_h-\sigma^{n}_h}{\tau}\right\|^2_{L^2}+C\tau^4+Ch^{2(r+1)},
\end{align}
where the second inequality can be obtained by the fact that
\begin{align*}
&\sigma^0_h=R_hu^0-U^0_h=0,\\
&{\|\sigma^1_h\|_{H^1}=\|R_hu^1-U^1_h\|_{H^1}=\|R_h(u^1-U^1)\|_{H^1}\qquad\quad\mbox{(use \eqref{FEMIN} and  \eqref{INITtime})}}\nonumber\\
&\qquad\,\,\,\,{\le C\|u^1-U^1\|_{H^1}\le C\tau^3.\qquad\qquad\qquad\quad\quad\,\,\mbox{(use \eqref{theory2} and \eqref{e1-error})}}
\end{align*}
Applying Gronwall's inequality and $\sigma^0_h=0$ to \eqref{sigmaresult}, there exists $\tau_*>0$, when $\tau\le\tau_*$
\begin{align}
\max_{0\le n\le N}\|\nabla \sigma^n_h\|_{L^2}\le C_6(\tau^2+h^{r+1}).
\end{align}

Due to $\sigma^n_h\in H^1_0$, we have
\begin{align}
\max_{0\le n\le N}\|\sigma^n_h\|_{L^2}\le C\max_{0\le n\le N}\|\nabla \sigma^n_h\|_{L^2}\le C_{7}(\tau^2+h^{r+1}).
\end{align}
Combining above result with \eqref{theory}, we have
\begin{align}
\max_{0\le n\le N}\|u^n-U^n_h\|_{L^2}\le& \max_{0\le n\le N}\|u^n-R_hu^n\|_{L^2}+\max_{0\le n\le m}\|\sigma^n_h\|_{L^2}\nonumber\\
\le&
Ch^{r+1}\max_{0\le n\le N}\|u^n\|_{H^{r+1}}+C_{7}(\tau^2+h^{r+1})\nonumber\\
\le&
C_{8}(\tau^2+h^{r+1}).
\end{align}
{Thus, let $h_*=h^*$ and $C^0>C_{8}$ in above result, the proof of Theorem \ref{theorem1} is complete.}\qed

\section{Numerical results}\label{sec6}

In this section, two numerical examples are presented to verify our theoretical analysis, where a free software (Freefem++) is used to perform all computations.

\begin{example}\label{ex1}
  Consider the following cubic nonlinear Schr\"{o}dinger equation with wave operator:
  \begin{equation}\label{example1-1}
    \left\{\begin{aligned}
      &u_{tt}-\Delta u+\mathbf{i}u_t+|u|^2u+w(\mathbf{x})u=g, &\mathbf{x}\in \Omega, 0<t\leq T,\\
      &u(\mathbf{x},t)=0,  &\mathbf{x}\in \partial\Omega,\\
      & u(\mathbf{x},0)=u_0(\mathbf{x}), \quad u_t(\mathbf{x},0)=u_1(\mathbf{x}), &\mathbf{x}\in \Omega,
    \end{aligned}\right.
  \end{equation}
\end{example}
where $w(x,y)=-x^2y^2$, and $\Omega = \{(x,y):(x-0.5)^2+(y-0.5)^2<0.5^2\}$. Moreover, the initial boundary conditions and the term $g$ in right-hand side are obtained by the
exact solution
\begin{equation}
  u(x,y) = 20e^{\mathbf{i}8t}(1+8t^2)x^2(1-x)y^2(1-y).
\end{equation}

A quasi--uniform triangulation is generated by FreeFEM++ with $M$ nodes distributed on the boundary of the circular domain $\Omega$.  In the numerical implementation, we use both the linear and quadratic finite element when doing the spatial discretization, the final time is set to be $T=1$.

\begin{figure}[!tbp]
\begin{minipage}[t]{0.5\linewidth}
\centering
\includegraphics[width=0.75\textwidth]{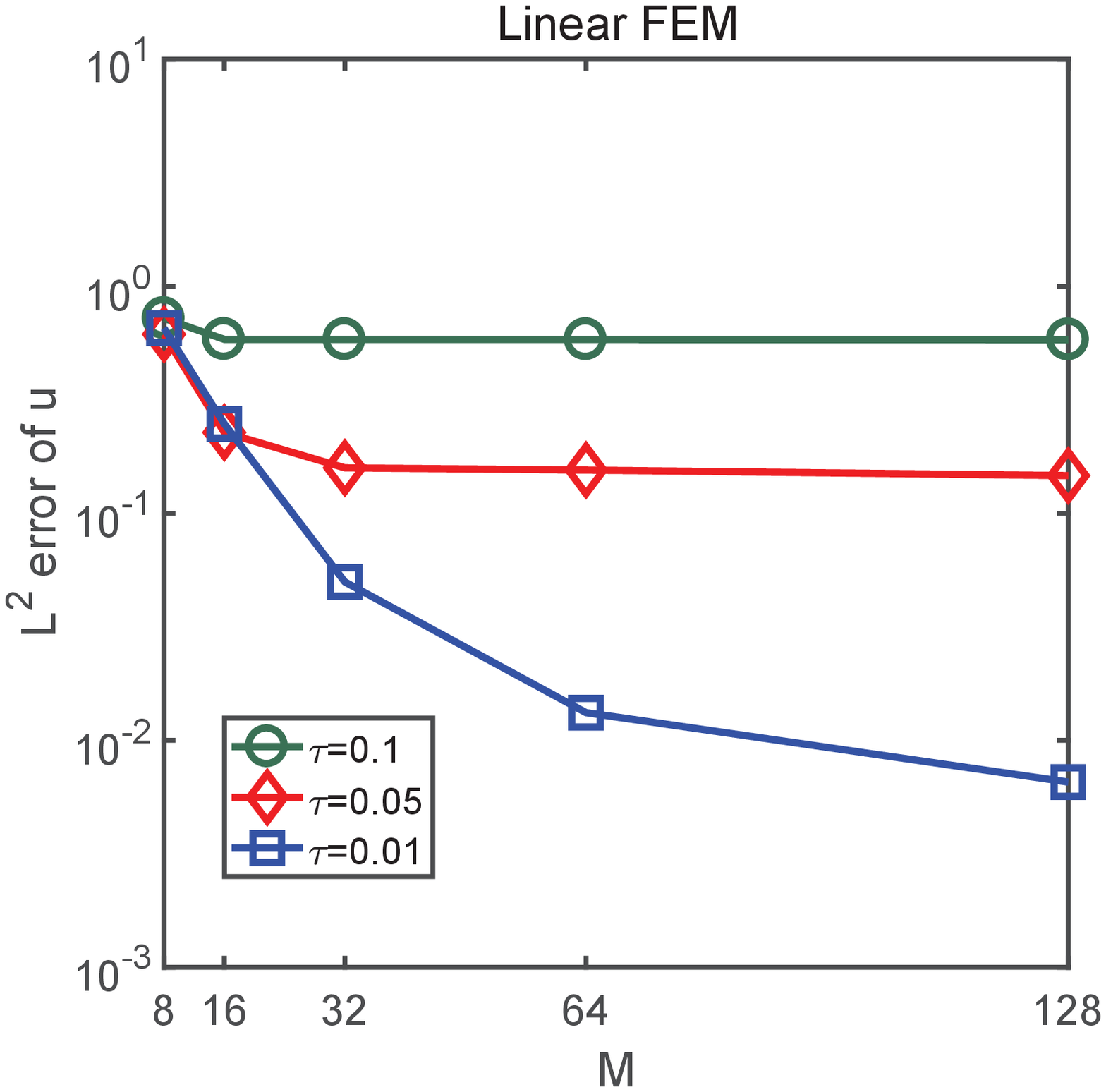}
\end{minipage}
\hfill
\begin{minipage}[t]{0.5\linewidth}
\centering
\includegraphics[width=0.75\textwidth]{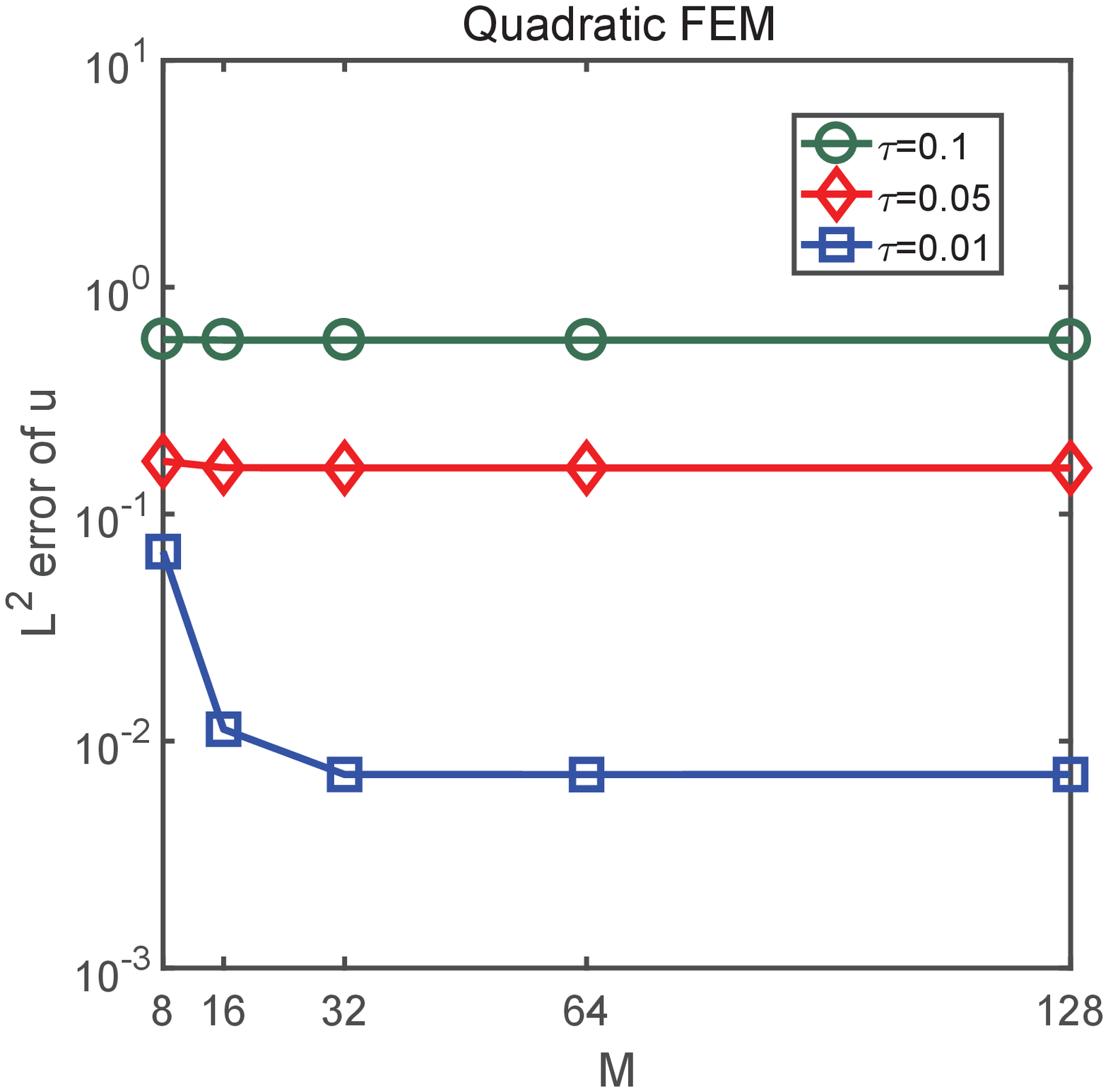}
\end{minipage}
\caption{\it{$L^2$--norm errors of the linear and quadratic FEMs for example 1.}}\label{fig1}
\end{figure}

We first list the $L^2$--norm errors of numerical solutions in Table \ref{table-1} with an extremely small time step size $\tau=2^{-14}$ so that the errors from temporal direction can be neglected. From Table \ref{table-1}, we can see that the numerical method reaches its optimal convergent order in spatial direction. Moreover, the convergence study for temporal direction using $L^2$--norm errors are also presented in Table \ref{table-2} with the small spatial step size $h=1/512$. From Table \ref{table-2}, we can see that the  numerical method  is second--order accurate in temporal direction. To check unconditional stability,  we take different spatial sizes $1/h=8, 16, 32, 64, 128$ with fixed time steps $\tau=0.1, 0.05, 0.01$. From Fig.~\ref{fig1}, we can see that, for each fixed time step size, both $L^2$--norm errors from the linear and quadratic FEMs  tend to be a constant when the mesh is refined gradually. It is shown that our numerical scheme is unconditionally stable. Thus, all numerical results completely match with our theoretical analysis.

\begin{table}[!tbp]
\centering
\caption{\small{Convergent analysis in spatial direction for example 1}}
\tabcolsep=40pt
\vspace{0.2cm}
\begin{tabular}{cccc}
\hline
$ $&P1 element  &P2 element \\
 \hline
$1/h=16$&2.509E$-$01&9.096E$-$03    \\
$1/h=32$&5.246E$-$02 &1.098E$-$03   \\
$1/h=64$&1.430E$-$02&1.486E$-$04      \\
\hline
$order_{avg}$&2.06 &2.97\\
\hline
\end{tabular}\label{table-1}
\end{table}

\begin{table}[!tbp]
\centering
\caption{\small{Convergent analysis in temporal direction for example 1.}}
\tabcolsep=40pt
\vspace{0.2cm}
\begin{tabular}{cccc}
\hline
$ $&P1 element  &P2 element \\
 \hline
$1/\tau=16$&2.456E$-$01&2.457E$-$01    \\
$1/\tau=32$&6.312E$-$02&6.322E$-$02   \\
$1/\tau=64$&1.598E$-$02&1.606E$-$03      \\
\hline
$order_{avg}$&1.97 &1.97\\
\hline
\end{tabular}\label{table-2}
\end{table}

\begin{example}\label{ex2}
  Consider the cubic nonlinear Schr\"{o}dinger equation with wave operator:
  \begin{equation}\label{example1-2}
    \left\{\begin{aligned}
      &u_{tt}-\Delta u+\mathbf{i}u_t+|u|^2u+w(\mathbf{x})u=0, &\mathbf{x}\in \Omega, 0<t\leq T,\\
      &u(\mathbf{x},t)=0,  &\mathbf{x}\in \partial\Omega,\\
      & u(\mathbf{x},0)=u_0(\mathbf{x}), \quad u_t(\mathbf{x},0)=u_1(\mathbf{x}), &\mathbf{x}\in \Omega,
    \end{aligned}\right.
  \end{equation}
\end{example}
where $w=-(\sqrt{2}\pi+\sin^2(\pi x) \sin^2(\pi y))$ and $\Omega =[0,1]\times[0,1]$. Moreover, the initial boundary conditions are obtained by the exact solution
\begin{equation}
  u(x,y) = \sin(\pi x)\sin(\pi y)e^{-{\bf i}\sqrt{2}\pi t}.
\end{equation}

Again, we solve the equation \eqref{example1-2} by both linear and quadratic FEMs up to time $T=1$. Similar to the above example, the convergent study in spatial direction using $L^2$--norm errors are shown in table \ref{table-3} by choosing extremely small time--step size $\tau=2^{-14}$.  Meanwhile, the convergent study in temporal direction  are presented with the small spatial step size $h=1/512$. Again, from Table \ref{table-3} and Table \ref{table-4}, we can see that the proposed method reaches its optimal convergent order in both  spatial and temporal directions.

\begin{figure}[!tbp]
\begin{minipage}[t]{0.5\linewidth}
\centering
\includegraphics[width=0.75\textwidth]{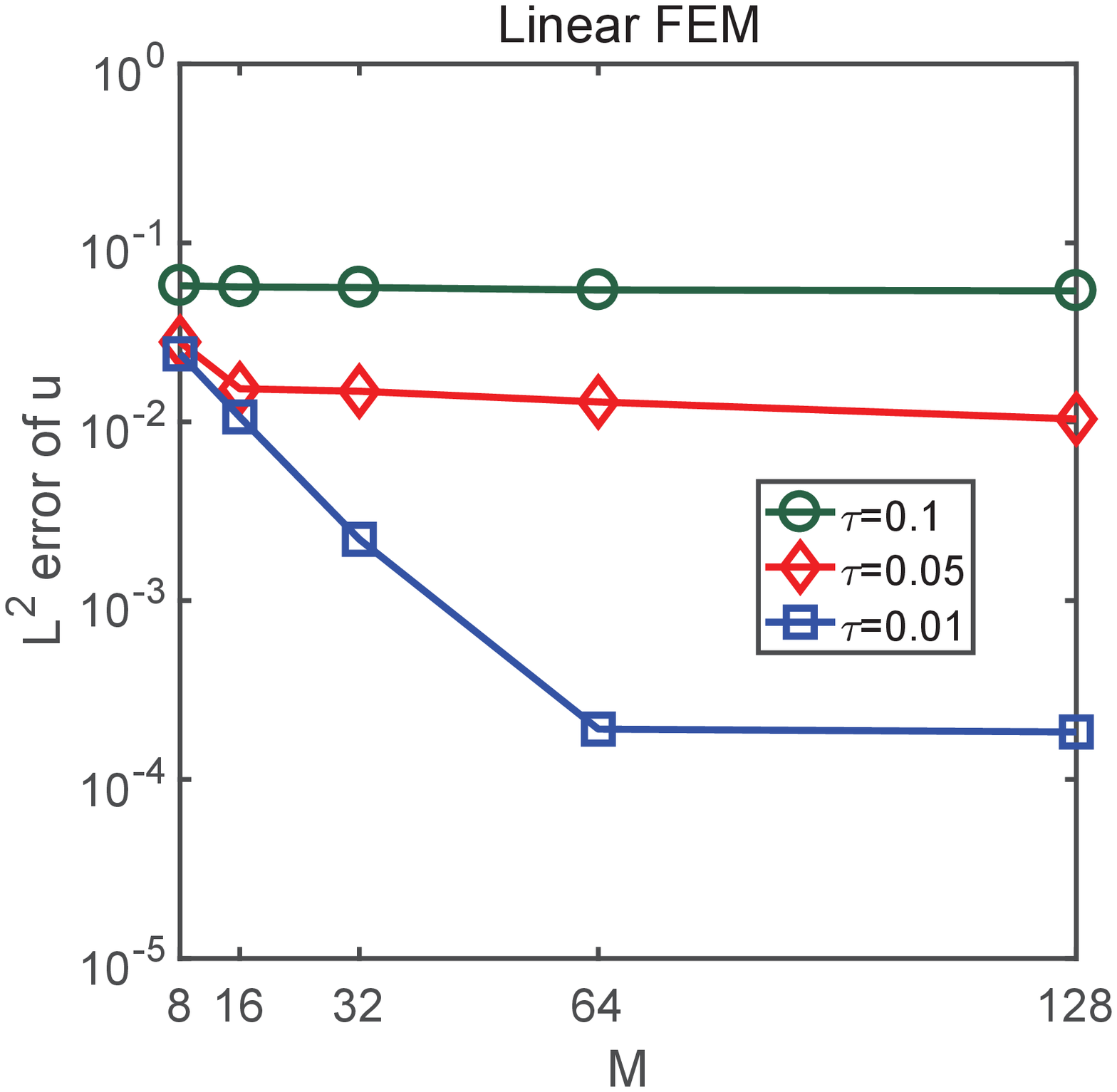}
\end{minipage}
\hfill
\begin{minipage}[t]{0.5\linewidth}
\centering
\includegraphics[width=0.75\textwidth]{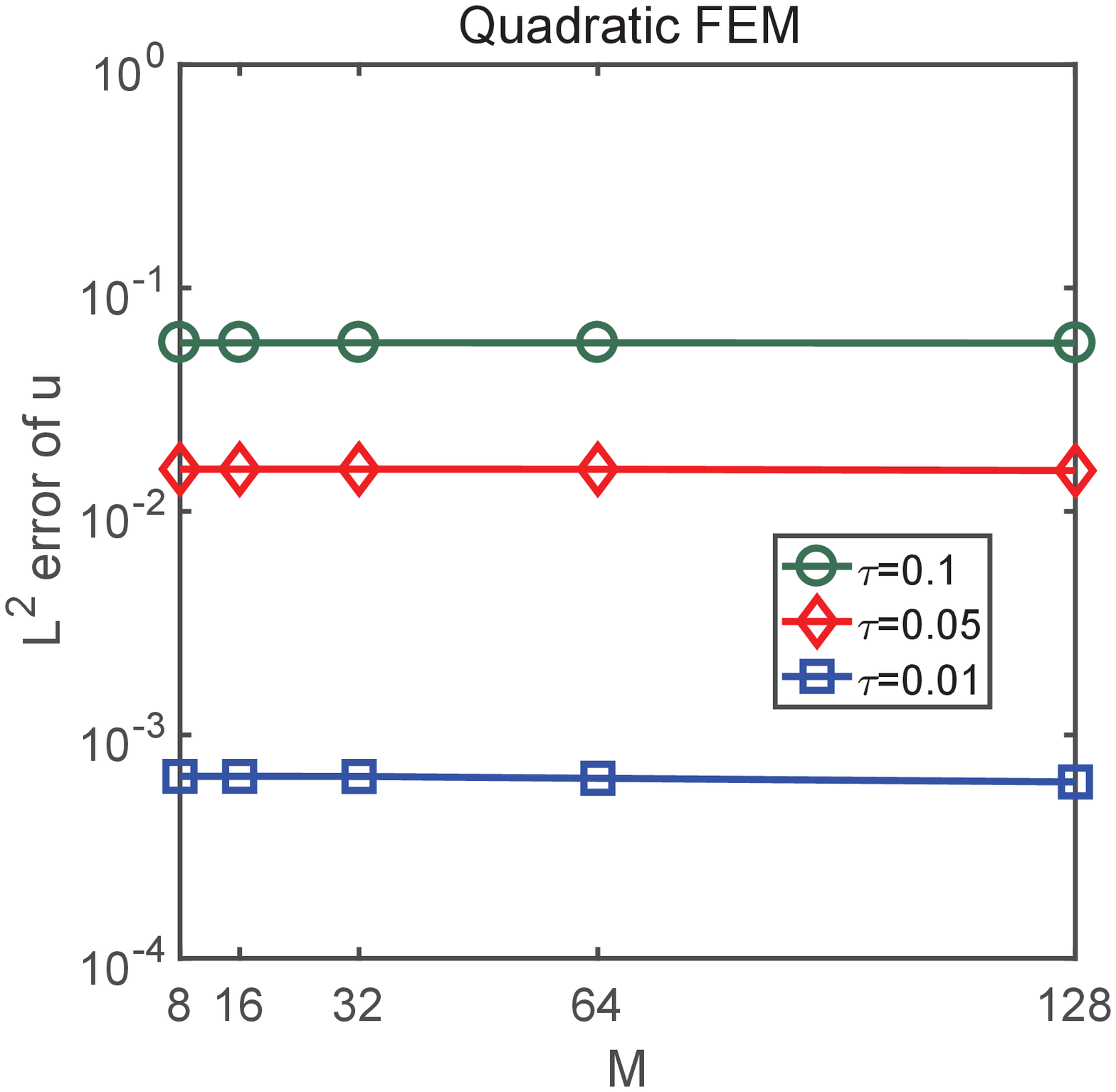}
\end{minipage}
\caption{\it{$L^2$--norm errors of the linear and quadratic FEMs for example 2.}}\label{fig2}
\end{figure}

Further, similar to the above example,  we gradually refine the spatial mesh size for three fixed time--step sizes $\tau=0.1, 0.05, 0.01$ at $T=1$, From Fig.~\ref{fig2}, we can see that the $L^2$--norm errors of both linear and quadratic FEMs tend to be a constant, which implies that our method is unconditionally stable.

\begin{table}[!tbp]
\centering
\caption{\small{Convergent analysis in spatial direction for example 2.}}
\tabcolsep=40pt
\vspace{0.2cm}
\begin{tabular}{cccc}
\hline
$ $&P1 element  &P2 element \\
 \hline
$1/h=16$&1.136E$-$02 &6.159E$-$05   \\
$1/h=32$&2.843E$-$03&7.511E$-$06      \\
$1/h=64$&7.126E$-$04&9.321E$-$07    \\
\hline
$order_{avg}$&2.00 &3.02\\
\hline
\end{tabular}\label{table-3}
\end{table}

\begin{table}[!tbp]
\centering
\caption{\small{Convergent analysis in temporal direction for example 2.}}
\tabcolsep=40pt
\vspace{0.2cm}
\begin{tabular}{cccc}
\hline
$ $&P1 element  &P2 element \\
 \hline
$1/\tau=16$&2.367E$-$02&2.365E$-$02    \\
$1/\tau=32$&6.196E$-$03&6.083E$-$03   \\
$1/\tau=64$&1.572E$-$03&1.503E$-$03      \\
\hline
$order_{avg}$&1.96 &1.99\\
\hline
\end{tabular}\label{table-4}
\end{table}

\begin{figure}[!tbp]
\centering
\includegraphics[width=0.5\textwidth]{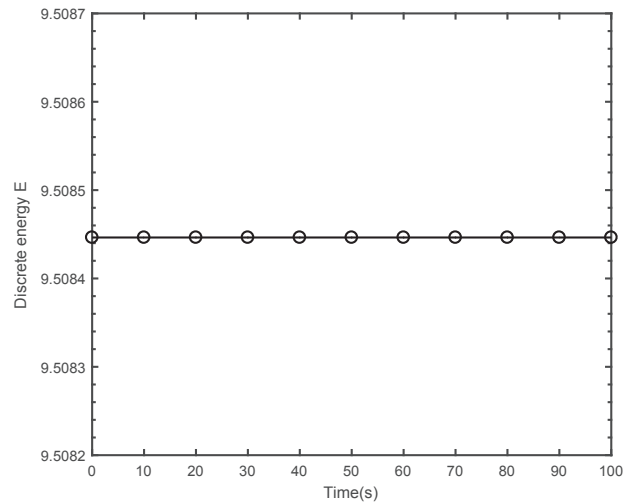}
\caption{\it{Evolution of discrete energy $E$ by linear finite element approximation for example 2.}}\label{fig3}
\end{figure}

Further, to illustrate the energy conservative property of our numerical scheme,  we present the numerical results at different time stages $T=0, 10, 20,\cdots,100,$ by using linear element approximation. From Fig.~\ref{fig3}, we can see that the discrete energy $E$ is conserved exactly (up to machine accuracy) with time evolution, which verifies our theoretical result in Theorem~\ref{Theorem1*}.

\section{Conclusions}
In this paper, an energy--conservative finite element method is present to solve nonlinear Schr\"{o}dinger equation with wave operator. Comparing to previous works \cite{CLC,CLC1}, our scheme is proved to keep energy conservation in a certain  discrete norm (Theorem \ref{Theorem1*}). Thus,  our scheme keeps higher stability. Moreover, we give unconditionally optimal error estimates of the modified leap--frog FEM for cubic nonlinear Schr\"{o}dinger equation with wave operator. 
By introducing a time--discrete system, the error of numerical solutions is split into two parts: the temporal error and the spatial error. We present the uniform boundedness of time--discrete solutions in some strong norms and the error estimates in temporal direction. Based on these results, we get the $L^2$ optimal error estimates in the sense that the time step size is not related to spatial mesh size. At last, numerical examples are provided for verifying the convergence--order, unconditional stability and  energy conservation of the proposed numerical method.

\section*{Acknowledgements}
Dongdong He was supported by  the president's fund--research start--up fund from the Chinese University of Hong Kong, Shenzhen (No. PF01000857). Kejia Pan was supported by the National Natural Science Foundation of China (Nos. 41874086, 41474103), the Excellent Youth Foundation of Hunan Province of China (No. 2018JJ1042) and the Innovation--Driven Project of Central South University (No. 2018CX042).

\end{document}